\documentclass{amsart}
\usepackage{amssymb}
\usepackage{amsfonts}
\usepackage{amsmath}
\usepackage{lastpage}
\usepackage[legalpaper,bookmarks=true,colorlinks=true,linkcolor=blue,citecolor=blue]{hyperref}
\usepackage{graphicx}%
\usepackage{fancyhdr}
\usepackage{color}
\usepackage[mathlines]{lineno}
\usepackage{lscape}
\usepackage{epsfig}
\usepackage{natbib}
\usepackage{geometry}
\usepackage{tgbonum}
\fontfamily{qcr}\selectfont

\newtheorem{theorem}{Theorem}
\theoremstyle{plain}

\newtheorem{corollary}{Corollary}

\newtheorem{lemma}{Lemma}

\numberwithin{equation}{section}

\newcommand{\Bin}{\bigskip \noindent}

\newcommand{\Ni}{\noindent}

\begin{document}
\title[Glivenko-Cantelli classes for arbitrary stationary variables]{Arbitrary functional Glivenko-Cantelli classes and applications to different types of dependence}

\author{Harouna Sangar\'e}
\author{Gane Samb Lo}
\author{Mamadou Cherif Moctar Traor\'e}

\begin{abstract}
Using a general strong law of large number proved by Sangar\'e and Lo (2015) and the entropy numbers, we provide functional Glivenko-Cantelli (GC) classes for arbitrary stationary real-valued random variables (rrv's). Next, the general results are particularized for different types of dependence (association, $\phi$-mixing, in particular) and compared with available results in the literature.

\bigskip
\noindent $^{\dag}$ Harouna Sangar\'e\\
DER MI,Facult\'e des Sciences et Techniques (FST), USTT-B, Mali (main affiliation)\\
LERSTAD, Gaston Berger University, Saint-Louis, S\'en\'egal.\newline
Institutional emails : harounasangare@fst-usttb-edu.ml, harouna.sangare@mesrs.ml, sangare.harouna@ugb.edu.sn\\
Own email : harounasangareusttb@gmail.com\\

\noindent $^{\dag\dag}$ Gane Samb Lo.\\
LERSTAD, Gaston Berger University, Saint-Louis, S\'en\'egal (main affiliation).\newline
LSTA, Pierre and Marie Curie University, Paris VI, France.\newline
AUST - African University of Sciences and Technology, Abuja, Nigeria\\
Institutional emails :gane-samb.lo@edu.ugb.sn, gslo@aust.edu.ng\\
Own email : ganesamblo@ganesamblo.net\\
Permanent address : 1178 Evanston Dr NW T3P 0J9, Calgary, Alberta, Canada.\\


\noindent \noindent $^{\dag\dag\dag}$ Mamadou Cherif Moctar Traor\'e.\\
LERSTAD, Gaston Berger University, Saint-Louis, S\'en\'egal.\newline
traore.cherif-mamadou-moctar@ugb.edu.sn.\\

\noindent\textbf{Keywords}. Association, $\phi$-mixing, Stationarity, Entropy number, Glivenko-Cantelli class\\
\textbf{AMS 2010 Mathematics Subject Classification:} 62G20; 60G10\\

\end{abstract}

\maketitle

\noindent

\section{Introduction} \label{sec01}

\noindent Let $\lbrace X,(X_n)_{n\geq 1}\rbrace$ be a sequence of real-valued random variables (rrv's) defined on the same probability space $(\Omega,\mathcal{A},\mathbb{P})$ and associated to the same cumulative distribution function (cdf) $F$. Let $\mathcal{L}_2(\mathbb{P}_X)$ be the class of all measurable functions $f: \mathbb{R}\rightarrow \mathbb{R}$ with $\mathbb{E}f(X)^2<+\infty$. Let us defined the functional empirical probability based on the $n$ first observations, $n\geq 1$, by

\begin{equation}\label{eq01}
\mathbb{P}_n(f)=\frac{1}{n}\sum_{i=1}^{n}f(X_i), \ f\in \mathcal{L}_2,
\end{equation}

\bigskip\noindent which will be centered at the mathematical expectations

$$
\mathbb{E}\mathbb{P}_n(f)=\mathbb{P}_X(f), \ f\in \mathcal{L}_2.
$$

\bigskip\noindent In this paper, we are interested by functional Glivenko-Cantelli classes, that is classes $\mathcal{F}\subset \mathcal{L}_2$ for which we have

$$
\sup_{f\in \mathcal{F}}|\mathbb{P}_n(f)-\mathbb{E}\mathbb{P}_n(f)|=\|\mathbb{P}_n-\mathbb{E}\mathbb{P}_n\|_{\mathcal{F}}\rightarrow 0 \ a.s.,\ as \ n\rightarrow +\infty.
$$

\bigskip\noindent We might also consider a class $\mathcal{C}$ of measurable sets in $\mathbb{R}$ and define

\begin{equation}\label{eq02}
\mathbb{P}_n(C)=\frac{1}{n}\sum_{i=1}^{n}Card\lbrace X_i\in C\rbrace, \ C\in \mathcal{C},
\end{equation}

\bigskip\noindent with

$$
\mathbb{E}\mathbb{P}_n(C)=\mathbb{P}_X(C)=\mathbb{P}(X\in C)
$$

\bigskip\noindent Taking $\mathcal{F}_c=\lbrace 1_{\rbrack-\infty,x\rbrack}, x\in \mathcal{R}\rbrace$ in Definition \eqref{eq01} or $\mathcal{C}_c=\lbrace \rbrack-\infty,x\rbrack, x\in \mathbb{R}\rbrace$ in Definition \eqref{eq02}, leads to the classical empirical function

$$
F_n(x)=\frac{1}{n}Card\lbrace i\in \lbrack 1,n\rbrack, X_i\leq x\rbrace, \ x\in \mathbb{R},
$$

\bigskip\noindent which in turn gives the Glivenko-Cantelli law for independent and identically distributed (\textit{iid}) sequences $(X_n)_{n\geq 1}$ under the form :

\begin{equation}\label{eq03}
\sup_{x\in \mathbb{R}}|F_n(x)-F(x)|\rightarrow 0\ a.s., \ as \ n\rightarrow +\infty.
\end{equation}

\bigskip\noindent Such a result, also known as the fundamental theorem of statistics is the frequentist paradigm (in opposition to the Bayesian paradigm) in statistics. In the form of \eqref{eq03}, the Glivenko-Cantelli law has gone trough a large number of studies for a variety of type of dependence. Also, it has been extensively used in statistical theory both for finding plug-in asymptotically efficient estimators and the related statistical tests based on the Donsker theorem that we will not study here.\\

\noindent To give a few examples, we cite the following results, \cite{billingsley} showed the convergence in law on $D\lbrack 0,1\rbrack$ of the empirical process for $\phi$-mixing rrv's under the condition $\sum_{k>0}k^2\sqrt{\phi(k)}<+\infty$. \cite{yosh} obtained the same result for $\alpha$-mixing under the condition on the mixing coefficient $\alpha_n=O(n^{-a})$ with $a>3$. This result has been first improved by \cite{shao} by only assuming $a>2$, then by \cite{shaoyu} who suppose that $a>1+\sqrt{2}$. \cite{rio} obtained a best condition $a>1$. Similar results are given by in \cite{shaoyu} for $\rho$-mixing rrv's, \cite{doukh} for $\beta$-mixing rrv's. For the associated dependence, \cite{yu93} obtained the convergence under the condition $Cov(X_1,X_n)=O(n^{-a})$, $a>7.5$. Next, \cite{shaoyu} weakened the covariance condition $a>\frac{3+\sqrt{33}}{2}$. \cite{louh} gave another proof and a result improvement with $a>4$.\\

\noindent As to the Glivenko-Cantelli type theorems for associated random variables, \cite{bagai} have first proposed an estimator for the survival function, discussed its asymptotic properties and then gave Glivenko-Cantelli theorem under moment condition, namely $\sum_{j=n+1}^{+\infty}\lbrack Cov(X_1,X_j)\rbrack^{1/3}<O(n^{-(r-1)})$, for some $r>1$, for stationary sequences in some compact subset. 
\cite{yu93} dropped the stationary assumption and considered a sequence of associated random variables having the same marginal distribution function. He obtained a Glivenko-Cantelli theorem under the following conditions : he first supposed that the distribution function is continuous and then

$$
\sum_{n=1}^{+\infty}\frac{1}{n^2}Cov(X_n,S_{n-1})<\infty
$$

\bigskip\noindent and if the sequence is stationary, then the last condition can be weakened to

$$
\frac{1}{n}\sum_{n=1}^{+\infty}Cov(X_n,X_i)\rightarrow 0, \ as \ n\rightarrow +\infty.
$$

\bigskip\noindent However, functional versions seem not to have been developed for dependent data while they are far more interesting than restricting to the particular case of $\mathcal{F}_c$. Besides, the functional version has the advantage to be linear in the sense that

$$
\forall (a,b)\in \mathbb{R}^2, \forall (f,g)\in (\mathcal{L}_1(\mathbb{P}_X))^2, \mathbb{P}_n(af+bg)=a\mathbb{P}_n(f)+b\mathbb{P}_n(g),
$$

\bigskip\noindent which allows the use of more mathematical latitudes.\\

\noindent In this paper, we directly address the functional Glivenko-Cantelli as introduced in Theorem in \cite{vdv}. We use general conditions of the covariances to establish GC-classes regardless the dependence structures. The obtained conditions are then applied to specific types of dependence.\\
 
\noindent The main achievements are Theorem \ref{hgc_theoApp}, as a very general functional Glivenko-Cantelli laws for arbitrary law under under Conditions \eqref{ghcf_01} and \eqref{ghcf_02} and Theorem \ref{thg-gc-ef-01} general Glivenko-cantelli laws for the real-valued empirical functions under Conditions \eqref{gcep1} and \eqref{gcep2}, also  to be check for different types of dependence. Next, we will focus on the second type of Glivenko-Cantelli classes for $\phi$-mixing and associated sequences with new results to be compared with previous results.\\

\noindent That functional approach needs the use of concentration numbers we define in the next section. We will also need to make a number of recalls on dependance types as $\phi$-mixing and associated sequences and other tools.\\

\noindent Taking that into account leads us to the following organization of the paper. The main result concerning the arbitrary stationary Glivenko-Cantelli classes will be stated in Section \ref{sec05}. But before that, Section \ref{sec02} will devote to entropy numbers and Vapnik-\v{C}ervonenkis classes, Section \ref{sec03} to recall of associated sequences, Section \ref{sec04} to the recall of the \cite{sanglo} general law which will be instrumental to our proofs. In Section \ref{sec06}, we focus on GC-classes regarding real-valued empirical function under types of dependence. \ref{sec07} concerns concluding remarks. The paper is ended by the Appendix, where is postponed the detailed proof of the main theorem, in Section \ref{sec08}. Reader familiar with entropy numbers or associated sequence may skip the related section to go directly to the section of their interest.\\

\section{Entropy numbers and Vapnik-\v{C}ervonenkis} \label{sec02}

\noindent It might be perceived that the following recall on the entropy numbers, which comes from combinatorial Theory, makes the paper heavier. But, in our view, it may help the reader who is not well aware of such techniques. Let $(E,\|\circ\|,\leq)$ be an ordered real normed space meaning that the order is compatible with the operations in the following sense

$$
\forall(\lambda,z)\in \mathbb{R}_+\vee\lbrace 0\rbrace\times E, \forall(x,y)\in E^2, (x<y)\rightarrow ((x+z<y+z) \ and \ (\lambda x<\lambda y)).
$$

\bigskip\noindent A bracket set of level $\varepsilon>0$ in $E$ is any set of the form

$$
B(\ell_1, \ell_2,\varepsilon)=\lbrace x\in E, \ell_1\leq x\leq \ell_2\rbrace, \ \ell_1\leq \ell_2\in E \ and \ \|\ell_2-\ell_1\|<\varepsilon.
$$

\bigskip\noindent For any subset $F$ of $E$, the bracketing entropy number at level $\varepsilon>0$, denoted $N_{(F,\|\circ\|,\varepsilon)}$ is the minimum of the numbers $p\geq 1$ for which we have $p$ bracket sets (or simply $p$ brackets) at level $\varepsilon$ covering $F$, where the $\ell_i$'s do not necessarily belong to $F$. The bracketing entropy number is closely related to the Vapnick-\v{C}ervonenkis index (VC-index) of a Vapnick-\v{C}ervonenkis set or class (VC-class or VC-set).\\

\noindent To define a VC-set, we need to recall some definitions. A subset $D$ of $B\subset E$ is picked out by a subclass $\mathcal{C}$ of the power set $\mathcal{P}(E)$ from $B$ if and only if $D$ is element of $\lbrace B\cap C, C\in \mathcal{C}\rbrace$. Next, $B$ is shattered by $\mathcal{C}$ if and only if all subsets of $B$ are picked out by $\mathcal{C}$ from $B$. Finally, the class $\mathcal{C}$ is a VC-class if and only if there exists an integer $n\geq 1$ such that no set of cardinality $n$ is shattered by $\mathcal{C}$. The minimum of those numbers $n$ minus one is the index of that VC-class. The VC-class is the most quick way to bound bracketing entropy numbers $N_{(F,\|\circ\|,\varepsilon)}$, as stated in \cite{vdv} : If $\mathcal{C}$ is a VC-class of index $I(\mathcal{C})$, then

\begin{equation}
N_{(F,\|\circ\|,\varepsilon)}\leq KI(\mathcal{C})(4e)^{I(\mathcal{C})}(1/\varepsilon)^{rI(\mathcal{C}-1)},
\end{equation}

\bigskip\noindent where $K>0$ and $r>1$ are universal constants.\\

\noindent Finally, for a class $\mathcal{F}$ of real-value functions $f : E\rightarrow \mathbb{R}$, we may defined the bracketing entropy number associated to $\mathcal{F}$ is the bracketing entropy number of the class $\mathcal{C}(\mathcal{F})$ of sub-graphs $f$ of elements $f\in \mathcal{F}$

$$
S_g(f)=\lbrace (x,t)\in E\times \mathbb{R}, f(x)>t\rbrace,
$$

\bigskip\noindent in $E^\ast=E\times \mathbb{R}$ endowed with the product norm which still is a normed space. We denoted by $N_{(F,\|\circ\|,\varepsilon)}$ to distinguish with the bracketing number using the norm of the functions $f\in \mathcal{F}$. As well, $\mathcal{F}$ is a VC-supg-class if and only if $\mathcal{C}(\mathcal{F})$ is a VC-class and $I(\mathcal{F})=I\mathcal{C}(\mathcal{F})$.\\

\noindent We will use such entropy numbers for formulating VC-classes.\\

\noindent Since the results are applied to associated sequences, we also proceed to a brief recall on them.\\

\section{Recall on associated sequences}\label{sec03}

\noindent To begin with, we remind some useful results on associated data that find in \cite{rao}.\\

\begin{lemma}\label{lem01}
(\cite{newm}). Suppose that $X,Y$ are two rv's with finite variance and, $f$ and $g$ are $\mathbb{C}^1$ complex valued functions on $\mathbb{R}^1$ with bounded derivatives $f^{\prime}$ and $g^{\prime}$. Then 

$$
|Cov(f(X),g(Y))|\leq \|f^{\prime}\|_{\infty}\|g^{\prime}\|_{\infty}Cov(X,Y).
$$
\end{lemma}

\Bin We also have :\\

\begin{lemma}\label{lem02}
Let $X$ and $Y$ be associated random variables with absolutely continuous distributions. Assume that the marginal densities $f_X$ and $f_Y$ are bounded by $M$. Then, for every $T>0$,

$$
H(x,y)=Cov(1_{\rbrack-\infty;x\rbrack}(X),1_{\rbrack-\infty;y\rbrack}(Y))\leq M^\ast\left(T^2Cov(X,Y)+\frac{1}{T}\right),
$$

\bigskip\noindent where $M^\ast=\max(\frac{2}{\pi^2},45M)$.\\
\end{lemma}

\bigskip\noindent Optimizing the choice of $T$ on the previous result, they find the following important inequality\\

\begin{corollary}
Under the same assumptions as in Lemma \ref{lem02} if $Cov(X,Y)>0$, one has that

$$
Cov(1_{\rbrack-\infty;x\rbrack}(X),1_{\rbrack-\infty;y\rbrack}(Y))\leq \frac{1}{M^\ast}Cov^{1/3}(X,Y).
$$
\end{corollary}

\begin{lemma} \label{bagaiprakasa}
(\cite{bagai}). Suppose the pair $X$ and $Y$ are associated random variables with bounded continuous densities $f_X$ and $f_Y$, respectively. Then there exists an absolute constant $c>0$ such that

$$
\sup_{x,y}|\mathbb{P}(X\leq x, Y\leq y)-\mathbb{P}(X\leq x)\mathbb{P}(Y\leq y)|\leq c\left\lbrace \max\left(\sup_{x}f_X(x),\sup_{x}f_Y(y)\right)\right\rbrace^{2/3}Cov^{1/3}(X,Y).
$$
\end{lemma}

\bigskip\noindent Finally, the main tool used in our main theorem is the Arbitrary Strong Law of Large number we describe below.

\section{The Sangar\'e-Lo SLNN, $\phi$-mixing and important other tools} \label{sec04}

\noindent We begin with this useful result proved by \cite{sanglo}. \\
\begin{lemma}(\cite{sanglo}). \label{lemmasanglo}
Let $X_1,X_2,...$ be an arbitrary sequence of rv's, and let $(f_{i,n})_{i\geq 1}$ be a sequence of measurable functions such that $\mathbb{V}ar\lbrack f_{i,n}(X_i)\rbrack<+\infty$, for $i\geq 1$ and $n\geq 1$. If for some $\delta$, $0<\delta<3$

\begin{equation}
C_1=\sup_{n\geq 1}\sup_{q\geq 1}\mathbb{V}ar\left(\frac{1}{q^{(3-\delta)/4}}\sum_{i=1}^{q}f_{i,n}(X_i)\right)<+\infty \label{ghc_01}
\end{equation}

\bigskip\noindent and for some $0<\delta<3$,

\begin{equation}
C_2=\sup_{n\geq 0}\sup_{k\geq 0}\sup_{q:q^2+1\geq k\leq(q+1)^2 }\sup_{k\geq j\leq(q+1)^2 }\mathbb{V}ar\left(\frac{1}{q^{(3-\delta)/2}}\sum_{i=1}^{j-q^2+1}f_{q^2+i,n}(X_{q^2+i})\right)<+\infty \label{ghc_02}
\end{equation}

\bigskip\noindent hold, then

$$
\frac{1}{n}\sum_{i=1}^{n}\left(f_{i,n}(X_i)-\mathbb{E}(f_{i,n}(X_i))\right)\rightarrow 0 \ a.s, \ as \ n\rightarrow +\infty.
$$
\end{lemma}

\bigskip \noindent \textbf{Remark.} We say that the sequence $X_{1}, X_{2},\cdots, X_{n}$ satisfies the \textbf{(GCIP)} whenever Conditions \eqref{ghc_01} and \eqref{ghc_02} hold, and we denote $\mathcal{C}$, the class of measurable functions for which \textbf{(GCIP)} holds.\\

\bigskip \noindent Since we are also going to apply our results to sequences verifying the $\phi$-mixing condition, we make a brief summary of this notion.

\subsection{$\phi$-mixing}\label{phimixing}

\noindent Let us have a brief recall on $\phi$-mixing. Let ($\Omega ,\mathcal{A},\mathbb{P}$) be a probability space and $\mathcal{A}_1$, $\mathcal{A}_2$ two sub $\sigma -$algebras of $\mathcal{A}$. The $\phi $-mixing coefficient is given by 

\begin{equation*}
\phi (\mathcal{A}_1,\mathcal{A}_2)=\sup \left\{ \left\vert \mathbb{P}(B)-\frac{\mathbb{P}(A\cap B)}{\mathbb{P}(A)}\right\vert, \
 A \in \mathcal{A}_1, \  B \in \mathcal{A}_2, \ \text{and}\ \mathbb{P}(A) \neq 0\right\} .
\end{equation*}

\noindent We remind that we have (see \cite{paul}) for $X$ $\mathcal{A}_1$-measurable and $Y$ $\mathcal{A}_2$-measurable, $q$-integrable and $p$-integrable respectively with 
$p>1$, $q>1$ and $(1/p)+(1/q)=1$, 

\begin{equation*}
Cov(X,Y)\leq 2\phi ^{1/p}(\mathcal{A}_1,\mathcal{A}_2)\left\vert \left\vert X\right\vert \right\vert
_{p}\left\vert \left\vert Y\right\vert \right\vert _{q}\text{ }p\text{, }
q\geq 1\text{ and }\frac{1}{p}+\frac{1}{q}
\end{equation*}

\Ni and

\begin{equation*}
0\leq \phi (\mathcal{A}_1,\mathcal{A}_2)\leq 1.
\end{equation*}

\noindent We define the strong mixing coefficient $\phi(n)$ by

\begin{equation*}
\phi (n)=\sup \left\{ \phi \left(\mathcal{F}_{1}^{k},\mathcal{F}_{n+k}^{\infty }\right), \ k\in \mathbb{N}^{\ast }\right\} 
\end{equation*}

\noindent where $\mathcal{F}_{j}^{\ell }$ is the $\sigma -$algebra generated by the variables $(X_{i}$, $j\leq i\leq \ell )$. We say that  $(X_{n})_{n\geq 1}$ is $\phi$-mixing if $\phi (n)\rightarrow 0$ as $n\rightarrow \infty$. For some further clarification on this point, the reader may have a quick look at \cite{paul}.\\

\noindent Let us just cite two useful analysis tools. \\

\subsection{Two useful lemmas} Next here is the Kronecker lemma. First, here is the Ces\`aro lemma.

\begin{lemma} \label{lemmacesaro} Let $(x_n)_{n\geq 1} \subset \mathbb{R}$ be a sequence of finite real numbers converging to $x \in \mathbb{R}$, then sequence of arithmetic means

$$
y_n=\frac{x_1+...+x_n}{n}, \ n\geq 1
$$

\noindent also converges to $x$.
\end{lemma}

\Bin For a proof, see also \cite{loeve} and \cite{ips-mfpt-ang} (page 367). Next, we recall the Kronecker lemma.

\begin{lemma} \label{lemmakronecker} (Kronecker Lemma). If $(b_n)_{n\geq 0}$ is an increasing sequence of positive numbers and $(x_n)_{n\geq 0}$ is a sequence of finite real numbers such that $\left( \sum_{1\leq k \leq n} x_k\right)_{n\geq 0}$ converges to a finite real number $s$, then 

$$
\frac{\sum_{1\leq k \leq n} b_k x_k}{b_n} \rightarrow 0 \ as \ n\rightarrow \infty. 
$$ 
\end{lemma}

\Bin For a proof, see \cite{loeve} or \cite{ips-mfpt-ang} (page 369).
 
\Bin In the next section, we provide general functional Glivenko-Cantelli classes.

\section{Our results} \label{sec05}

\noindent In this section, we give general functional Glivenko-Cantelli classes. These results will be used in Section \ref{sec06} to establish GC-classes for real-valued empirical function under types of dependence.

\subsection{General functional GC-classes}

\noindent Let us begin by giving a slight different version of Theorem 2.4.1 in \cite{vdv}, page 122.

\begin{theorem}\label{hgc_theo} Let $X, X_{1}, X_{2}, \cdots$ be an arbitrary stationary sequence of rrv's with common \textit{cdf} $F$.\\
	
\noindent a) Let $\mathcal{C}$ be a class of measurable set such that \\
	
\noindent (a1) For any $\varepsilon > 0$, $N_{[]}\left(\mathcal{C}, \left| \circ \right|, \varepsilon \right) < + \infty$,\\
	
\noindent (a2) For any $C \in \mathcal{C}$, $\mathbb{P}_{n} \left(C\right) \longrightarrow \mathbb{P}_{X}\left(C\right)$, as $n \longrightarrow +\infty$.\\

\bigskip \noindent Then $\mathcal{C}$ is GC-class, that is
	
\begin{equation}\label{hgc_02}
\lim_{n \to +\infty} \sup_{C \in \mathcal{C}}\left|\mathbb{P}_{n}\left(C\right) - \mathbb{P}_{X} \left( C \right) \right|=0.
\end{equation}
	
\bigskip \noindent b) Let $\mathcal{F}$ be a class of measurable set such that\\
	
\noindent (b1) For any $\varepsilon > 0$, $N_{[]} \left( \mathcal{F}, \| \circ \|_{\mathcal{L}_{2}\left(\mathbb{P}_{X} \right)}, \varepsilon \right) < +\infty$,\\
	
\noindent (b2) For any $f \in \mathcal{F}$, $\mathbb{P}_{n}\left(f\right) \longrightarrow \mathbb{P}_{X}\left(f\right)$, as $n \to + \infty$.\\
	
\bigskip \noindent Then $\mathcal{F}$ is GC-class, that is
	
\begin{equation}\label{hgc_03}
\lim_{n \to +\infty} \sup_{f \in \mathcal{F}}\left| \mathbb{P}_{n} \left( f\right) - \mathbb{P}_{X} \left( f \right) \right|=0.
\end{equation}
\end{theorem}

\Bin \textbf{Proof}. The proof of the Theorem \ref{hgc_theo} is postponed in the Appendix (Section \ref{sec08}) $\square$.\\

\noindent \textbf{Remark (R2)}. If we apply Part (a) to $\mathcal{C}_c$, Condition (a1) can be dropped since it is VC-class of index 2. Indeed, for 
$A=\{x_1, \ x_2, x_3\}$ with $x<x_2<x_3$. The subset $\{x_2\}$ cannot be picked out by $\mathcal{C}_c$ from $A$ since, for any $x\in \mathbb{R}$,  
$]-\infty, x]\cap A$ will be of on the five sets : $\emptyset$ (for $x\leq x_1$), $\{x_1\}$ (for $x_1<x\leq x_2$), $\{x_1, x_2\}$ (for $x_2<x\leq x_3$) and $A$ (for $x>x_3$).\\

\noindent Then finding GC-classes for $\mathcal{C}_c$ reduces to establishing Condition (b1). For now, we are going to focus on the application of GC-classes of the functional empirical process for $\mathcal{C}_c$, that is, on results as in Formula \eqref{eq03}. The focus will be on types of dependence, given it is known that $\mathcal{C}_c$ is GC-class for \textit{iid} data. Using the \textit{GCIP} conditions in Lemma \ref{lemmasanglo} leads to our applicable results as follows.

\begin{theorem}\label{hgc_theoApp} Let $X, X_{1}, X_{2}, \cdots$ be an arbitrary stationary sequence of rrv's with common \textit{cdf} $F$.\\
	
\noindent a) Let $\mathcal{C}$ be a class of measurable set such that \\
	
\noindent (a1) For any $\varepsilon > 0$, $N_{[]}\left(\mathcal{C}, \left| \circ \right|, \varepsilon \right) < + \infty$,\\
	
\noindent (a2) For any $C \in \mathcal{C}$, the following conditions hold : for some $0<\delta<3$,
	
\begin{equation}
\sup_{q\geq 1}\dfrac{1}{q^{(3-\delta )/2}}\sum_{i,j=1}^{q}\left[ 
\mathbb{P(}X_{i}\in C,X_{j}\in C)-\mathbb{P(}X_{i}\in C)\mathbb{P}(
X_{j}\in C)\right] <+\infty. \label{ghcC_01}
\end{equation}

\bigskip\noindent and 

\begin{equation}
\sup_{k\geq 1}\sup_{q\text{ : }q^{2}+1\leq k\leq
(q+1)^{2}}\sup_{k\leq j\leq (q+1)^{2}}\dfrac{1}{q^{(3-\delta )}}%
\sum_{i,\ell=1}^{j-q^{2}+1}\left[ 
\begin{array}{c}
\mathbb{P(}X_{_{q^{2}+i}}\in C,X_{_{q^{2}+\ell}}\in C) \\ 
-\mathbb{P(}X_{_{q^{2}+i}}\in C)\mathbb{P(}X_{_{q^{2}+\ell}}\in C)%
\end{array}
\right] <+\infty . \label{ghcC_02}
\end{equation}
\bigskip \noindent Then $\mathcal{C}$ is GC-class, that is
	
\begin{equation}\label{hgc_02}
\lim_{n \to +\infty} \sup_{C \in \mathcal{C}}\left|\mathbb{P}_{n}\left(C\right) - \mathbb{P}_{X} \left( C \right) \right|=0.
\end{equation}
	
\bigskip \noindent b) Let $\mathcal{F}$ be a class of measurable set such that\\
	
\noindent (b1) For any $\varepsilon > 0$, $N_{[]} \left( \mathcal{F}, \| \circ \|_{\mathcal{L}_{2}\left(\mathbb{P}_{X} \right)}, \varepsilon \right) < +\infty$,\\
	
\noindent (b2) For any $f \in \mathcal{F}$, the following conditions hold : for some $0<\delta<3$,

\begin{equation}
\sup_{q\geq 1}\mathbb{V}ar\left(\frac{1}{q^{(3-\delta)/4}}\sum_{i=1}^{q}f(X_i)\right)<+\infty \label{ghcf_01} 
\end{equation}

\bigskip\noindent and 

\begin{equation}
\sup_{k\geq 0}\sup_{q:q^2+1\geq k\leq(q+1)^2 }\sup_{k\geq j\leq(q+1)^2 }\mathbb{V}ar\left(\frac{1}{q^{(3-\delta)/2}}\sum_{i=1}^{j-q^2+1}f(X_{q^2+i})\right)<+\infty. \label{ghcf_02}
\end{equation}

\bigskip \noindent Then $\mathcal{F}$ is GC-class, that is
	
\begin{equation}\label{hgc_03}
\lim_{n \to +\infty} \sup_{f \in \mathcal{F}}\left| \mathbb{P}_{n} \left( f\right) - \mathbb{P}_{X} \left( f \right) \right|=0.
\end{equation}
\end{theorem}

\Bin \textbf{Proof}. It follows the lines of that of Theorem \ref{hgc_theo} in the Appendix (Section \ref{sec08})

\subsection{GC-classes for the empirical function for dependent data}

\noindent As already mentioned, let us focus on $\mathcal{F}_c$, on 

$$
F_{n}(x)=\frac{1}{n}\sum_{i=1}^{n} 1_{]-\infty, \ x]}(X_i), \ x \in \mathbb{R},
$$

\Ni and

$$
F(x)=\mathbb{E}\left(1_{]-\infty, \ x]}(X_i)\right), \ x \in \mathbb{R}.
$$

\Bin Now, we apply Part (a) of Theorem \ref{hgc_theoApp} for $C_x=]-\infty, x]$, $x \in \mathbb{R}$ to have

\begin{theorem}  \label{thg-gc-ef-01} Let $X$, $X_{1}$ ,$X_{2}$, $\cdots$,  be an arbitrary stationary and square integrable sequence of \textit{rrv}'s, defined on a probability space $(\Omega, \mathcal{A}, \mathbb{P})$. We have 

$$
\sup_{x\in \mathbb{R}}|F_n(x)-F(x)|\rightarrow 0 \ a.s., \ as \ n\rightarrow +\infty
$$

\Bin whenever the following general conditions hold : for  some $\delta \in ]0;3[$

\begin{equation}
\sup_{q\geq 1}\dfrac{1}{q^{(3-\delta )/2}}\sum_{i,j=1}^{q}\left[ 
\mathbb{P(}X_{i}\leq x,X_{j}\leq x)-\mathbb{P(}X_{i}\leq x)\mathbb{P}(
X_{j}\leq x)\right] <+\infty. \label{gcep1}
\end{equation}
and
\begin{equation} \label{gcep2}
\sup_{k\geq 1}\sup_{q\text{ : }q^{2}+1\leq k\leq
(q+1)^{2}}\sup_{k\leq j\leq (q+1)^{2}}\dfrac{1}{q^{(3-\delta )}}%
\sum_{i,\ell=1}^{j-q^{2}+1}\left[ 
\begin{array}{c}
\mathbb{P(}X_{_{q^{2}+i}}\leq x,X_{_{q^{2}+\ell}}\leq x) \\ 
-\mathbb{P(}X_{_{q^{2}+i}}\leq x)\mathbb{P(}X_{_{q^{2}+\ell}}\leq x)%
\end{array}
\right] <+\infty .
\end{equation}
\end{theorem}

\noindent \textbf{Remark (R3)}. Since the expressions in \eqref{gcep1} and \eqref{gcep2} are zero whenever the variables are independent and identically distributed (iid), the iid case holds without any further condition, which is the classical Glivenko-Cantelli theorem.

\Bin \textbf{Proof}. We have to apply  Part (a) of Theorem \ref{hgc_theoApp} for $C_x=]-\infty, x]$, $x \in \mathbb{R}$. The conditions become

\begin{eqnarray*}
&&\sup_{q\geq 1}\mathbb{V}ar\left( \dfrac{1}{q^{(3-\delta )/4}}%
\sum_{i=1}^{q}1_{\left\{ X_{i}\leq x\right\} }\right) =\sup_{q\geq 1}\dfrac{1}{q^{(3-\delta )/2}}\mathbb{V}ar\left( \sum_{i=1}^{q}1_{\left\{ X_{i}\leq
x\right\} }\right)  \notag\\
&&=\sup_{q\geq 1}\dfrac{1}{q^{(3-\delta )/2}}\sum_{i,j=1}^{q}\left[ 
\mathbb{P(}X_{i}\leq x,X_{j}\leq x)-\mathbb{P(}X_{i}\leq x)\mathbb{P}(X_{j}\leq x)\right] <+\infty.
\end{eqnarray*}

\Bin and

\begin{eqnarray*}
&&\sup_{k\geq 1}\sup_{q\text{ : }q^{2}+1\leq k\leq (q+1)^{2}}\sup_{k\leq
j\leq (q+1)^{2}}\mathbb{V}ar\left( \dfrac{1}{q^{(3-\delta )/2}}%
\sum_{i=1}^{j-q^{2}+1}1_{\left\{ X_{_{q^{2}+i}}\leq x\right\} }\right) 
\notag \\
&=&\sup_{k\geq 1}\sup_{q\text{ : }q^{2}+1\leq k\leq (q+1)^{2}}\sup_{k\leq
j\leq (q+1)^{2}}\dfrac{1}{q^{(3-\delta )}}\mathbb{V}ar\left(
\sum_{i=1}^{j-q^{2}+1}1_{\left\{ X_{_{q^{2}+i}}\leq x\right\} }\right) 
\notag \\
&=&\sup_{k\geq 1}\sup_{q\text{ : }q^{2}+1\leq k\leq
(q+1)^{2}}\sup_{k\leq j\leq (q+1)^{2}}\dfrac{1}{q^{(3-\delta )}}%
\sum_{i,\ell=1}^{j-q^{2}+1}\left[ 
\begin{array}{c}
\mathbb{P(}X_{_{q^{2}+i}}\leq x,X_{_{q^{2}+\ell}}\leq x) \\ 
-\mathbb{P(}X_{_{q^{2}+i}}\leq x)\mathbb{P(}X_{_{q^{2}+\ell}}\leq x)%
\end{array}%
\right] <+\infty 
\end{eqnarray*}

$\square$.\\

\Ni Let us focus on two types of dependence : association and $\phi$-mixing.

\section{Application to the real-valued empirical function} \label{sec06}

\noindent In the present section, we use the Theorem \ref{thg-gc-ef-01} to give applications to association and $\phi$-mixing. 

\subsection{Associated case} \label{subsec_sec_01_sec_06}

\begin{corollary} \label{theo-assoc_01}
Suppose that the $X_j$'s are associated and form a stationary and square integrable and have bounded continuous \textit{pdf}, then  

$$
\sup_{x\in \mathbb{R}}|F_n(x)-F(x)|\rightarrow 0 \ a.s., \ as \ n\rightarrow +\infty
$$

\Bin holds whenever we have

\begin{equation}
\lim_{q\rightarrow +\infty }\frac{1}{q}\sum_{j=2}^{q}Cov^{1/3}(X_{1},X_{j})\rightarrow 0. \label{gcep2-assoc}
\end{equation}
\end{corollary}

\Bin \textbf{Proof}. $X_j$'s are associated (see Section (\ref{sec03}, page\pageref{sec03})). By applying Lemma (\ref{bagaiprakasa}, page \pageref{bagaiprakasa}), conditions \eqref{gcep1} and \eqref{gcep2} become the following one, for $\nu =(1-\delta)/2\geq 0$ with $0<\delta <1$, 

\begin{equation}
\sup_{q\geq 1}\dfrac{1}{q^{1+\nu }}C\sum_{i,j=1}^{q}Cov^{1/3}(X_{i},X_{j})<+%
\infty \label{gc-a}
\end{equation}

\noindent and

\begin{equation}
\sup_{q\geq 1}\dfrac{1}{q^{2\left( 1+\nu \right) }}C\sum_{i,j=q^{2}+1}^{\left( q+1\right) ^{2}}Cov^{1/3}(X_{i},X_{j})<+\infty . \label{gc-b}
\end{equation}

\Ni Now since the sequence is second order stationary, then \eqref{gc-a} and \eqref{gc-b} will be equivalent to

\begin{equation*}
\sup_{q\geq 1}\dfrac{1}{q^{\nu }}\left[ \mathbb{V}ar\left( 1_{\left\{
X_{i}\leq x\right\} }\right) +\frac{2}{q}%
\sum_{j=2}^{q}(q-j+1)Cov^{1/3}(X_{1},X_{j})\right] <+\infty.
\end{equation*}

\noindent Which in turn becomes, by the kronecker lemma (see Lemma \ref{lemmakronecker}, page \pageref{lemmakronecker}),

\begin{equation*}
\mathbb{V}ar\left( 1_{\left\{ X_{i}\leq x\right\} }\right) +\frac{2}{q}%
\sum_{j=2}^{q}(q-j+1)Cov^{1/3}(X_{1},X_{j})<+\infty
\end{equation*}

\noindent and finally

\begin{equation*}
\lim_{q\rightarrow +\infty }\frac{1}{q}\sum_{j=2}^{q}Cov^{1/3}(X_{1},X_{j})
\rightarrow 0
\end{equation*}
$\square$.

\Bin We rediscover the result by \cite{yu93}.

\begin{corollary} \label{theo-assoc_01}
Suppose that the $X_j$'s are associated and form a stationary and square integrable and have bounded continuous \textit{pdf}, then  

$$
\sup_{x\in \mathbb{R}}|F_n(x)-F(x)|\rightarrow 0 \ a.s., \ as \ n\rightarrow +\infty
$$

\Bin holds whenever we have

\begin{equation}
\lim_{q\rightarrow +\infty }\frac{1}{q}\sum_{j=2}^{q}Cov(X_{1},X_{j})\rightarrow 0. \label{gcep2-assoc}
\end{equation}
\end{corollary}

\Bin \textbf{Proof}. Since $F$ is strictly increasing and continuous, we can use the sequence $U_{i}=F(X_{i})$, $i\geq 1$. The  $U_{i}$'s are obviously  associated, stationary and square integrable $(0,1)$-uniformly distributed random variables on $[0, \ 1]$ and we have for each $n\geq 1$

\begin{equation*}
\forall x\in \mathbb{R}, \ F_{n}(x)=G_{n}(F(x)),
\end{equation*}

\noindent where $G_{n}(x))$ is the real empirical function based on $U_{1}$, $U_{2}$, $\cdots$, $U_{n}$. So

$$
\sup_{x\in \mathbb{R}}|F_n(x)-F(x)|\rightarrow 0 \ a.s., \ as \ n\rightarrow +\infty,
$$

\Bin holds if and only if

\begin{equation}
\sup_{q\geq 1}\frac{1}{q^{1+\nu }}\mathbb{V}ar\left(
\sum_{i=1}^{q}U_{i}\right) <+\infty  \label{gca1}
\end{equation}

\Bin and

\begin{equation}
\sup_{q\geq 1}\frac{1}{q^{2(1+\nu )}}\mathbb{V}ar\left(
\sum_{i=q^{2}+1}^{(q+1)^{2}}U_{i}\right) <+\infty ,  \label{gcal2}
\end{equation}

\noindent where $\nu =(1-\delta )/2\geq 0$ with $0<\delta <1$. That are \begin{equation}
\sup_{q\geq 1}\frac{1}{q^{1+\nu }}\mathbb{V}ar\left(
\sum_{i=1}^{q}F(X_{i})\right) <+\infty  \label{gc1}
\end{equation}

\noindent and

\begin{equation}
\sup_{q\geq 1}\frac{1}{q^{2(1+\nu )}}\mathbb{V}ar\left(
\sum_{i=q^{2}+1}^{(q+1)^{2}}F(X_{i})\right) <+\infty .  \label{gc2}
\end{equation}

\noindent If the sequence is second order stationary, then (\ref{gc1}) implies (\ref{gc2}), since

\begin{eqnarray*}
\frac{1}{q^{2(1+\nu )}}\mathbb{V}ar\left(
\sum_{i=q^{2}+1}^{(q+1)^{2}}F(X_{i})\right) &=&\frac{\left( 2q+1\right)
^{1+\nu }}{q^{2(1+\nu )}}\left[ \frac{1}{\left( 2q+1\right) ^{1+\nu }}%
\mathbb{V}ar\left( \sum_{i=1}^{2q+1}F(X_{i})\right) \right] \\
&\sim &\frac{2}{q^{(1+\nu )}}\mathbb{V}ar\left( \frac{1}{k^{\left( 1+\nu
\right) /2}}\sum_{i=1}^{k}F(X_{i})\right) ,
\end{eqnarray*}

\noindent for $k=2q+1.$ And (\ref{gc1}) may be written as
 
\begin{equation}
\sup_{q\geq 1}\frac{1}{q^{\nu }}\left[ \mathbb{V}ar(F(X_{1}))+\frac{2}{q}%
\sum_{i=2}^{q}\left( q-i+1\right) Cov(F(X_{1}),\text{ }F(X_{i}))\right]
<+\infty .  \label{q2}
\end{equation}

\noindent This is our general condition under which we have the Glivenko-Cantelli class for second order stationary associated sequence. Then, by the
Kronecker lemma (see Lemma \ref{lemmakronecker}, page \pageref{lemmakronecker}), we have the Glivenko-Cantelli class if 

\begin{equation}
\sigma ^{2}=\mathbb{V}ar(F(X_{1}))+2\sum_{i=2}^{+\infty }Cov(F(X_{1}),\text{ 
}F(X_{i}))<+\infty .  \label{gc01}
\end{equation}

\noindent Clearly, by the Ces\`{a}ro lemma (see Lemma \ref{lemmacesaro}, page \pageref{lemmacesaro}), Formula (\ref{gc01}) implies

\begin{equation}
\lim_{q\rightarrow +\infty }\frac{1}{q}\sum_{i=1}^{q}Cov(F(X_{1}),\text{ }%
F(X_{i}))\rightarrow 0.  \label{gc02}
\end{equation}

\noindent By the Lemma (\ref{lem01}, page \pageref{lem01}), since $\sup_{x}\left\vert f_{X}(x)\right\vert $ is
finite, (\ref{gc02}) reduces to

\begin{equation*}
\lim_{q\rightarrow +\infty }\frac{1}{q}\sum_{i=1}^{q}Cov(X_{1}, \ X_{i})\rightarrow 0
\end{equation*}
$\square$.\\

\noindent Such a result is also obtained by \cite{yu93} for the strong convergence of empirical distribution function for associated
sequence with identical and continuous distribution.

\subsection{$\phi $-mixing case}

\noindent We already made a brief recall on $\phi$-mixing in Subsection (\ref{phimixing}) in Section (\ref{sec04}). We are now going to provide applications to it.\\

\begin{corollary}
Suppose that the $X_j$'s form a $\phi$-mixing  stationary and square integrable sequence of random variables with mixing condition $\phi$. Then

$$
\sup_{x\in \mathbb{R}}|F_n(x)-F(x)|\rightarrow 0 \ a.s., \ as \ n\rightarrow +\infty.
$$

\Bin holds whenever we have

\begin{equation}
\phi (r)=O\left( r^{\frac{-4}{1-\delta }}\right)  \label{gcep2-phi}
\end{equation}
\end{corollary}

\noindent \textbf{Proof}. In this case, our conditions become

\begin{eqnarray*}
&&\sup_{q\geq 1}\dfrac{1}{q^{(3-\delta )/2}}\mathbb{V}ar\left(
\sum_{i=1}^{q}1_{\left\{ X_{i}\leq x\right\} }\right)  \\
&=&\sup_{q\geq 1}\dfrac{1}{q^{(3-\delta )/2}}\left(\sum_{i=1}^{q}\mathbb{V}%
ar\left( 1_{\left\{ X_{i}\leq x\right\} }\right) +2\sum_{1\leq i<j\leq
q}Cov\left( 1_{\left\{ X_{i}\leq x\right\} },1_{\left\{ X_{j}\leq x\right\}
}\right) \right) \\
&\leq &\dfrac{q}{q^{(3-\delta )/2}}+\dfrac{2}{q^{(3-\delta )/2}}\sum_{1\leq
i\leq q-1;j=i+h,h<r}\phi (\left\vert i-j\right\vert ) \\
&&+\dfrac{2}{q^{(3-\delta )/2}}\sum_{1\leq i\leq q-1;j=i+h,h\geq r}\phi
(\left\vert i-j\right\vert )F_{i}(x)F_{j}(x) \\
&\leq &\frac{q+2(q-1)(r-1)}{q^{(3-\delta )/2}}+\frac{q^{2}-q-2(q-1)(r-1)}{%
q^{(3-\delta )/2}}\phi (r) \\
&\leq &\frac{q(2+2r-1)-2r+2+q^{2}\phi (r)}{q^{(3-\delta )/2}} \\
&\sim &\frac{2qr}{q^{(3-\delta )/2}}+\frac{q^{2}\phi (r)}{q^{(3-\delta )/2}}%
<+\infty .
\end{eqnarray*}

\noindent Finally, our condition reduces to

\begin{equation*}
r=\left[ q^{\frac{1-\delta }{2}}\right] \text{ and }\phi (r)=O\left(
q^{-2}\right) =O\left( r^{\frac{-4}{1-\delta }}\right)
\end{equation*}
$\square$.

\Bin \textbf{Important remark}. In the case of $\phi$-mixing sequences, we do not need whole function $\phi$. Instead we may fix $x \in \mathbb{R}$ and consider the modulus $\phi(x,r)$ related to the sequence $1_{]-\infty,x]}(X_i)$'s. It is clear that 

$$
\forall r\geq 1, \ \forall x \in \mathbb{R}, \ \phi(x,r)\leq \phi(r) 
$$

\Bin We have the Glivenko-Cantelli theorem if and only if

$$
\forall x \in \mathbb{R}, \ \phi(x,r)=O\left( r^{\frac{-4}{1-\delta }}\right).
$$

\Bin The same can be done for any particular class $\mathcal{C}$ by using the associated sequences $1_{C}(X_i)$'s and the $\phi(C,r)$ mixing modulus, and get the condition

$$
\forall x \in \mathcal{C} , \ \phi(C,r)=O\left( r^{\frac{-4}{1-\delta }}\right).
$$

\section{Conclusion} \label{sec07}
\noindent In a coming paper, we will focus on the functional Glivenko-Cantelli classes with a considerable classes of functions and its applications to a number of situations.

\section{Appendix} \label{sec08}

\noindent\textbf{Proof of the Theorem \ref{hgc_theo}.} Consider a class $\mathcal{F}$ of measurable real functions such that each $f\in \mathcal{F}$.
Let $p=N_{[]}(\mathcal{F},L_{1},\epsilon )$ be finite, for every $\epsilon >0$. Then there is a sequence of intervals $[f_{i},g_{i}]$, $i=1,\cdots,p$
such that for every $f_{i},\ g_{i}\in \mathcal{F}$, $1\leq i\leq p$, and $\mathcal{F}\subseteq \bigcup [f_{i},g_{i}],$ and such that $\mathbb{EP}%
_{n}(g_{i})-\mathbb{EP}_{n}(f_{i})\leq \epsilon $. For any $f\in \mathcal{F}$, there exist $f_{i}$ and $g_{i}$ such as $f_{i}\leq f\leq g_{i}$ and 

\begin{equation}
\mathbb{P}_{n}(f)-\mathbb{EP}_{n}(f)=\mathbb{P}_{n}(f)-\mathbb{EP}%
_{n}(g_{i})+\mathbb{EP}_{n}(g_{i})-\mathbb{EP}_{n}(f).
\end{equation}

\Bin On the other hand we have, by the monotonicity of the probability, that 

\begin{equation*}
\mathbb{P}_{n}(f_{i})\leq \mathbb{P}_{n}(f)\leq \mathbb{P}_{n}(g_{i}).
\end{equation*}

\Bin Therefore 

\begin{eqnarray*}
\mathbb{P}_{n}(f)-\mathbb{EP}_{n}(f) &\leq &\mathbb{P}_{n}(g_{i})-\mathbb{EP}%
_{n}(g_{i})+\mathbb{EP}_{n}(g_{i})-\mathbb{EP}_{n}(f) \\
&\leq &(\mathbb{P}_{n}-\mathbb{EP}_{n})(g_{i})+\epsilon .
\end{eqnarray*}

\Bin Next 

\begin{eqnarray*}
\mathbb{EP}_{n}(f)-\mathbb{P}_{n}(f) &\leq &\mathbb{EP}_{n}(f)-\mathbb{EP}%
_{n}(f_{i})+\mathbb{EP}_{n}(f_{i})-\mathbb{P}_{n}(f) \\
&\leq &\epsilon +(\mathbb{EP}_{n}-\mathbb{P}_{n})(f_{i}).
\end{eqnarray*}

\Bin Thus 

\begin{eqnarray}
||\mathbb{EP}_{n}-\mathbb{P}_{n}||_{\mathcal{F}} &=&\sup_{f\in \mathcal{F}}|%
\mathbb{EP}_{n}-\mathbb{P}_{n}|(f)  \label{p} \\
&=&\sup_{f\in \mathcal{F}}\max ((\mathbb{EP}_{n}-\mathbb{P}_{n})(f),(\mathbb{%
P}_{n}-\mathbb{EP}_{n})(f))  \notag \\
&\leq &\epsilon +\max_{1\leq i\leq p}\max ((\mathbb{EP}_{n}-\mathbb{P}%
_{n})(f_{i}),(\mathbb{P}_{n}-\mathbb{EP}_{n})(g_{i})).  \notag
\end{eqnarray}

\Bin Since for every $i$, $1\leq i\leq p$, we have $\mathbb{P}_{n}(f_{i})%
\longrightarrow \mathbb{EP}_{n}(f_{i})$ by the Lemma (\ref{lemmasanglo}, page \pageref{lemmasanglo}), it follows
that 

\begin{equation*}
\max_{1\leq i\leq p}\max ((\mathbb{EP}_{n}-\mathbb{P}_{n})(f_{i}),(\mathbb{P}%
_{n}-\mathbb{EP}_{n})(g_{i}))\longrightarrow 0.
\end{equation*}

\Bin By applying the upper limit on the left-hand side in (\ref{p}), we get for
every $\epsilon >0$ 

\begin{equation*}
\limsup_{n\rightarrow \infty }||\mathbb{EP}_{n}-\mathbb{P}_{n}||_{\mathcal{F}%
}\leq \epsilon ,\ a.s.
\end{equation*}
$\square$.

\end{document}